\documentclass[11pt]{amsart}
\usepackage{epsfig,amsmath}

\newtheorem{thm}[subsection]{Theorem}
\newtheorem{lem}[subsection]{Lemma}

\newtheorem{obs}[subsection]{Observation}

\theoremstyle{definition}
\newtheorem{definition}[subsection]{Definition}

\newcommand{\R}{\mathbf{R}}
\newcommand{\Z}{\mathbf{Z}}

\def\L{\Lambda}
\def\eps{\epsilon}
\def\x{\times}
\def\E{\mathcal{E}}
\def\Q{\mathcal{Q}}
\def\DA{\textit{DA}}
\def\ra{\rightarrow}

\def\ds{\displaystyle}

\begin{document}

\parskip 6pt
\parindent 0pt
\baselineskip 14pt

\title[No skew branes on quadrics\dots]{No skew
branes on non-degenerate hyperquadrics}

\author[Sha]{Ji-Ping Sha}
\author[Solomon]{Bruce Solomon}
\address{Department of Mathematics, Indiana University,
Bloomington, IN 47405}
\email{jsha@indiana.edu, solomon@indiana.edu}

\date{First draft November 12, 2004. Last Typeset \today.}

\begin{abstract}
We show that non-degenerate hyperquadrics in $\,\R^{n+2}\,$
admit no skew branes. Stated more traditionally, a compact
codimension-one immersed submanifold of a non-degenerate
hyperquadric of euclidean space must have parallel tangent
spaces at two distinct points. Similar results have been
proven by others, but (except for ellipsoids in $\,\R^3\,$)
always under $\,C^2\,$ smoothness and genericity
assumptions. We use neither assumption here.
\end{abstract}

\maketitle

\section{Introduction and preliminaries}

A \emph{non-degenerate hyperquadric} in $\,\R^{n+2}\,$ is a
level set of a {non-degenerate} real quadratic form on
$\,\R^{n+2}\,$. We show here that any $\,C^1\,$ immersion
from a compact smooth manifold $\,M^n\,$ into $\,\R^{n+2}\,$
with image on such a hyperquadric must have \emph{a pair of
parallel tangent planes}. Our main result expresses this
precisely:

\begin{thm}\label{thm:main}
Let $\,M\,$ denote a compact smooth $n$-manifold. If a
$\,C^1\,$ immersion $\,F:M\to\R^{n+2}\,$ maps $\,M\,$ into a
non-degenerate hyperquadric, then there are distinct points
$\,p,q\in M\,$ such that $\,dF(T_pM)=dF(T_qM)\,$.
\end{thm}

In short, the unoriented gauss map on a compact immersed
$n$-dimensional submanifold of a non-degenerate hyperquadric
in $\,\R^{n+2}\,$ cannot be injective. \cite{gs} and \cite{t}
call a compact immersed codimension two submanifold of
$\,\R^{n+2}\,$ with injective unoriented gauss map a
\emph{skew brane}. From that viewpoint, Theorem
\ref{thm:main} says:

\begin{quotation}\it
{No skew brane in
$\,\R^{n+2}\,$ lies on a non-degenerate hyperquadric}.
\end{quotation}

Compare this fact with the main result of \cite{gs}: Namely
that in $\,\R^3\,$, 1-dimensional skew branes---\emph{skew
loops}---can be immersed on any closed $\,C^2\,$ surface
which is \emph{not} quadric. It seems likely that skew branes
are common in higher odd dimensions too.

Other authors have obtained weaker versions of Theorem
\ref{thm:main}. In 1971, generalizing an earlier
curve-theoretic result of B.~Segre \cite{se}, J.H.~White
ruled out ``generic'' $\,C^2\,$ skew branes on ellipsoids
\cite{w}. In 2002, S.~Tabachnikov extended White's result to
all other non-degenerate hyperquadrics under the same
genericity and smoothness hypotheses \cite{t}. More recently,
M.~Ghomi \cite{g} was able to remove White's $\,C^2\,$ and
genericity assumptions, but only for the case of
1-dimensional skew branes---\emph{skew loops}---on ellipsoids
in $\,\R^3\,$.

Our main contribution here lies with Lemmas \ref{obs:reg} and
\ref{lem:retract}, which allow us to remove the genericity
and $\,C^2\,$ assumptions quite generally. We thank our
colleague Paul Kirk for asking a question which suggested the
latter Lemma. That result helps us improve an earlier argument
by the second author, which proved Theorem \ref{thm:main} for
all non-degenerate hyperquadrics except ellipsoids.

We now discuss some preparatory ideas and notation.

Given a non-degenerate quadric hypersurface $\,\Q\,$ in
$\,R^{n+2}\,$, one can always find a real symmetric
non-singular matrix $\,Q\,$ and a number $\,c\,$ such that
\[
\Q=\left\{x\in\R^{n+2}\colon x\cdot Q\, x = c\right\}\ .
\]
If $\,c\ne 0\,$, $\,\Q\,$ forms a smooth hypersurface in
$\,\R^{n+2}\,$. If $\,c=0\,$, $\,\Q\,$ forms a cone over a
smooth product of ellipsoids, and has an isolated singular
point at the origin. To distinguish the two cases, we shall
respectively say that $\,\Q\,$ is \emph{smooth} or
\emph{conical}.

Fixing $\,\Q\,$ and $\,Q\,$ as above, let $\,M\,$ be a compact
smooth $n$-manifold , and  $\,F:M\to\R^{n+2}\,$ a $\,C^1\,$
immersion whose image lies on $\,\Q\,$. We consider the energy
function $\,\E_F:M\x M\to\R\,$ defined by
\begin{equation*}
\E_F\left(p,q\right):=
\big(F(p)-F(q)\big)\cdot Q\, \big(F(p)-F(q)\big).
\end{equation*}
Following \cite{w} and \cite{t}, we can relate the existence
of parallel tangent planes for $\,F\,$ to that of certain
critical points for $\,\E_F\,$:

\begin{obs}\label{obs:key}
If $\,d\E_F\,$ vanishes at $\,(p,q)\in M\x M\,$, then either
$\,F\,$ has parallel tangent planes at $\,p\,$ and $\,q\,$,
or else $\,F(p)=\lambda\,F(q)\,$ for some $\,\lambda\in\R\,$,
with $\,\lambda = \pm 1\,$ if $\,\Q\,$ is smooth.
\end{obs}

\begin{proof}
Differentiating the identity $\,F\cdot Q\, F\equiv c\,$, one
sees immediately that $\,dF(T_pM)\,$ is perpendicular to
$\,Q\,F(p)\,$ and $\,dF(T_qM)\,$ is perpendicular to
$\,Q\,F(q)\,$.

At the same time, by differentiating $\,\E_F\,$, we similarly
find that both $\,dF(T_pM)\,$ and $\,dF(T_qM)\,$ are
perpendicular to $\,Q\,F(p)-Q\,F(q)\,$ whenever $\,d\E_F\,$
vanishes at $\,(p,q)\,$. Together, these facts make both
$\,dF(T_pM)\,$ and $\,dF(T_qM)\,$ perpendicular to both
$\,Q\,F(p)\,$ and $\,Q\,F(q)\,$ at any critical point of
$\,\E_F\,$.

Since $\,dF(T_pM)\,$ and $\,dF(T_qM)\,$ both have codimension
two in $\,\R^{n+2}\,$, we now get $\,dF(T_pM)=dF(T_qM)\,$,
i.e., a pair of parallel tangent planes, as long as
$\,Q\,F(p)\,$ and $\,Q\,F(q)\,$ are linearly independent.

If, on the other hand, they are linearly dependent,
then $\,Q\,F(p)=\lambda\,Q\,F(q)\,$ for some
$\,\lambda\in\R\,$. Since $\,Q\,$ is non-singular, this means
$\,F(p)=\lambda\,F(q)\,$ as well. And when $\,\Q\,$ is smooth,
we can further deduce that $\,\lambda=\pm 1\,$ by dotting the
preceding identities with each other and using $\,F(p)\cdot
Q\,F(p) = F(q)\cdot Q\,F(q)=c\ne 0\,$.
\end{proof}

When we prove our main result in \S\ref{sec:proof} below, the
preceding Observation will, virtually by itself, give the
desired conclusion when $\,\Q\,$ is conical. The smooth case,
however, requires the additional ideas we develop next.

Define three subsets of the product $\,M\x M\,$, two of which
depend on the immersion $\,F:M\to \Q\,$, as follows:
\begin{eqnarray*}
\Delta\,\, &:=& \left\{(p,p)\colon p\in M\right\}\quad\text{(diagonal of $\,M\x M\,$)}\ ,\\
D_F &:=& \left\{(p,q)\colon p\ne q\,,\ F(p)=F(q)\right\}\quad\text{(double-point locus)}\ ,\\
A_F &:=& \left\{(p,q)\colon F(p) = -F(q)\right\}\qquad\text{(antipodal locus)}
\end{eqnarray*}
The diagonal $\,\Delta\,$ is clearly a compact, embedded
submanifold of $\,M\x M\,$, and is diffeomorphic to $\,M\,$.

\begin{obs}\label{obs:disjoint}
When $\,\Q\,$ is smooth, $\,D_F\,$ and $\,A_F\,$ are disjoint
compact subsets of $\,(M\x M)\setminus\Delta\,$.
\end{obs}

\begin{proof}
We have $\,0\not\in F(M)\,$ because $\,\Q\,$ is smooth, and
hence $\,D_F\cap A_F=\emptyset\,$. Since $\,M\x M\,$ is
compact, the Observation follows easily if both $\,D_F\,$ and
$\,A_F\,$ are closed in $\,M\x M\,$. That $\,A_F\,$ is closed
is obvious. To see that $\,D_F\,$ is also closed, note that
$\,D_F\cup \Delta\,$ is closed; it then suffices to show
$\,\Delta\,$ is open in $\,D_F\cup\Delta\,$, which holds
because $\,F\,$ is an immersion. Indeed, each point $\,p\in
M\,$ has an open neighborhood $\,U_p\,$ on which $\,F\,$ is
injective, and hence $\,D_F\cap (U_p\x U_p)=\emptyset\,$.
\end{proof}

Unlike the diagonal $\,\Delta\,$, the double-point and
antipodal loci $\,D_F\,$ and $\,A_F\,$ are not, in general,
submanifolds of $\,M\x M\,$. We therefore make the following
definition:

\begin{definition}
We call a $\,C^1\,$ immersion $\,F:M\to\Q\,$
\emph{DA-regular} if $\,D_F\,$ and $\,A_F\,$ are both
embedded $(n-1)$-dimensional $\,C^1\,$ submanifolds of
$\,M\x M\,$.
\end{definition}

\section{Proof of the main theorem.}\label{sec:proof}

We precede our proof of Theorem \ref{thm:main} with two
Lemmas. The first was perhaps overlooked in
\cite{w} and \cite{t}. It will free us from the need to
\emph{assume} \DA -regularity of $\,F\,$ in Theorem
\ref{thm:main}. The second is a $\,C^1$-Morse-theoretic
statement that lets us dispense with the $\,C^2\,$ assumption
used by earlier authors.

\begin{lem}\label{obs:reg}
When $\,\Q\,$ is smooth, a $\,C^1\,$ immersion $\,F:M\to\Q\,$
is either DA-regular, or else has a pair of parallel
tangent planes.
\end{lem}

\begin{proof}
We will use the following general fact: Suppose
$\,T_1\,$ and $\,T_2\,$ are subspaces of a vector space
$\,T\,$. Then for any $\,(X,Y)\in T_1\x T_2\,$, we have
\begin{equation}\label{eqn:decomp}
\left(X,Y\right) = \textstyle{1\over
2}\,\left(X+Y,\,X+Y\right) + \textstyle{1\over
2}\,\left(X-Y,\,Y-X\right)\ .
\end{equation}
This corresponds to the direct sum decomposition
$\,T\x T=T_\Delta\oplus T_\Delta^\perp\,$, where
$\,T_\Delta:=\left\{(Z,Z)\colon Z\in T\right\}\,$,
and $\,T_\Delta^\perp:= \{(Z,-Z)\colon Z\in T\}\,$.
Noticing that the map $\,Z\mapsto (Z,-Z)\,$ gives an isomorphism
$\,T\approx T_\Delta^\perp\,$, we see that
$\,\{(X-Y,Y-X)\colon
X\in T_1,Y\in T_2\}=T_\Delta^\perp\,$ iff $\,T_1+
T_2 = T\,$, i.e. $\,T_1\,$ and $\,T_2\,$ are transverse in
$\,T\,$. It then follows easily from (\ref{eqn:decomp}) that
$\,T_1\x T_2\,$ is transverse to $\,T_\Delta\,$ in $\,T\x
T\,$ iff $\,T_1\,$ and $\,T_2\,$ are transverse in $\,T\,$.

To exploit this fact, consider the $\,C^1\,$ map
\begin{equation*}
\widehat F:(M\x M)\setminus\Delta\to\Q\x\Q\ ,\qquad \widehat
F\left(p,q\right):= \big(F(p),F(q)\big)\ .
\end{equation*}
Suppose that whenever $\,(p,q)\in D_F\,$, $\,dF(T_pM)\ne
dF(T_qM)\,$. As unequal hyperplanes, these two subspaces
are transverse in $\,T_{F(p)}\Q\,$. Taking
$\,T_1=dF(T_pM)\,$, $\,T_2=dF(T_qM)\,$ and $\,T=T_{F(p)}\Q =
T_{F(q)}\Q\,$ for each $\,(p,q)\in D_F\,$ now yields, via the
general fact above, transversality of $\,\widehat F\,$ with
respect to the diagonal $\,\Q_\Delta:=\left\{(x,x)\colon
x\in\Q\right\}\,$. This makes $\,D_F = \widehat
F^{-1}\left(\Q_\Delta\right)\,$ an embedded submanifold of
$\,M\x M\,$ with codimension $\,\dim M+\dim M -\dim\Q = 2n-(n+1)
= n-1\,$, as desired (\cite[Theorem 1.3.3]{h}).

To prove the same for $\,A_F\,$, note that $\,d(-F)(T_pM) =
dF(T_pM)\,$ for all $\,p\in M\,$, and run the same argument
with $\,\widehat F\,$ replaced by the map $\,(p,q)\mapsto
\big(F(p),-F(q)\big)\,$.
\end{proof}
\medskip

\begin{lem}\label{lem:retract}
Let $\,W\,$ be a connected compact  $\,C^\infty\,$ manifold, and let
$\,f:W\to[0,1]\,$  be a $\,C^1\,$ function. Suppose
$\,\L_0:=f^{-1}(0)\,$ and $\,\L_1:=f^{-1}(1)\,$ are
both non-empty $\,C^1$-embedded submanifolds of $\,W\,$. If
$\,f\,$ has no critical point on $\,W\setminus
(\L_0\cup\L_1)\,$, then $\,\L_1\,$ is a deformation retract
of $\,W\setminus \L_0\,$.
\end{lem}

\begin{proof}
Fix a smooth riemannian metric on $\,W\,$.
When an embedded  submanifold $\,\L\subset W\,$ is merely
$\,C^1\,$, the ``nearest-point retraction'' may not be
well-defined on any neighborhood of $\,\L\,$. Still, by
using a $\,C^1\,$ approximation of the $\,C^0\,$ normal plane
distribution along $\,\L\,$, one can nevertheless construct a \emph{tubular
neighborhood} of $\,\L\,$ (\cite[4.5]{h}). In
particular, we can find neighborhoods $\,U_0\,$ and $\,U_1\,$
of the submanifolds $\,\L_0\,$ and $\,\L_1\,$ respectively,
such that $\,U_0\cap U_1 = \emptyset\,$, and
\begin{itemize}
\item[(a)] $\,\L_1\,$ is a deformation retract of $\,U_1\,$,
while $\,W\setminus U_0\,$ is a deformation retract of $\,W\setminus\L_0\,$.
\end{itemize}

Because $\,C^\infty(W)\,$ is dense in $\,C^1(W)\,$
(\cite[2.2]{h}), we may now select an $\,\eps>0\,$ and a function $\,\tilde f\in
C^\infty(W)\,$ which approximates $\,f\,$ well enough to
ensure the following:
\begin{itemize}
\item[(b)]
$\ds{W\setminus U_0\subset\tilde f^{-1}\left([\eps,1+\eps]\right)}\,$,
$\ds{\tilde f^{-1}\left([1-\eps,1+\eps]\right)\subset U_1}\,$, and
\medskip

\item[(c)]
The set $\,\tilde f^{-1}\left([\eps,1-\eps]\right)\,$ contains no
critical point of $\,\tilde f\,$.
\end{itemize}

By (c) and standard Morse theory (\cite[Theorem 3.1]{mi}), we
deduce that
\begin{itemize}
\item[(d)]
$\,\tilde f^{-1}\left([1-\eps,1+\eps]\right)\,$ is a
deformation retract of $\,\tilde
f^{-1}\left([\eps,1+\eps]\right)\,$\ .
\end{itemize}

We then obtain the desired deformation retraction
$\,W\setminus\L_0\to\L_1\,$ with ease: First deform
$\,W\setminus\L_0\,$ onto $\,W\setminus U_0\,$ by (a); then
continue to deform it onto $\,\tilde
f^{-1}\left([1-\eps,1+\eps]\right)\,$ by (b) and (d); finally
deform onto $\,\L_1\,$ by (a) and (b).
\end{proof}

We are now ready to prove our main result.

\begin{proof}[Proof of Theorem \ref{thm:main}]
Consider first the more elementary case of conical $\,\Q\,$
where a quick calculation yields
\begin{equation}\label{eqn:epq}
\E_F(p,q) = -2F(p)\cdot Q\,F(q)\qquad\text{($\Q\,$
conical)}\ .
\end{equation}
To get a contradiction, suppose $\,F\,$ has no pair of
parallel tangent planes. Then at any critical point
$\,(p,q)\,$ of $\,\E_F\,$, Observation \ref{obs:key} gives
$\,F(p)=\lambda\,F(q)\,$ for some $\,\lambda\in \R\,$. Hence
$\,F(p)\cdot Q\,F(q) = \lambda F(q)\cdot Q\,F(q) = 0\,$. But
then $\,\E_F(p,q)=0\,$ by (\ref{eqn:epq}) above. This makes
$\,0\,$ the \emph{only} critical value of $\,\E_F\,$, so that
$\,\E_F\,$ must vanish identically on $\,M\x M\,$. But then
$\,d\E_F\,$ vanishes identically too, and Observation
\ref{obs:key} now says that $\,F(p)=\lambda F(q)\,$ for
\emph{all} $\,p,q\in M\,$. In particular, this confines
$\,F(M)\,$ to the line spanned by $\,F(q)\,$ for any $\,q\in
M\,$, which is clearly impossible. Our Theorem therefore holds
in the conical case.

To get the result for smooth $\,\Q\,$, note first that we can
assume connectivity of $\,M\,$. We will again seek a
contradiction by supposing $\,F\,$ has no parallel tangent
planes. By Lemma \ref{obs:reg} $\,A_F\,$ is a
embedded submanifold of codimension
$\,n+1\ge 2\,$ in the connected manifold $\,M\x M\,$, so that
$\,(M\x M)\setminus A_F\,$ remains connected.

Now consider the function $\,f:={1\over 4c}(1-\E_F)\,$, which has
the same critical points as $\,\E_F\,$. By Observation
\ref{obs:key}, the absence of parallel tangent planes forces all
critical points of $\,f\,$ to lie on $\,A_F=f^{-1}(0),\,$ and
$\,\Delta\cup D_F=f^{-1}(1)\,$. In particular, $\,f\,$ is
extremized on precisely these sets, and Lemma \ref{lem:retract}
now applies, making $\,\Delta\cup D_F\,$ a deformation retract
of $\,(M\x M)\setminus A_F\,$. The connectivity of the latter
now implies that of the former, and by Observation
\ref{obs:disjoint} this forces $\,D_F=\emptyset\,$.

On the other hand, applying Lemma \ref{lem:retract} again, this
time with $\,f:={1\over 4c}\E_F\,$, we may deduce that
$\,A_F\,$ is a deformation retract of $\,(M\x M)\setminus
\Delta\,$, and the isomorphisms of cohomology groups (using
$\,\Z_2\,$ coefficients henceforth)
\begin{equation*}
H^n\big((M\x M)\setminus\Delta\big)\approx H^n(A_F)\approx 0\
,
\end{equation*}
because $\,\dim A_F=n-1\,$.
Consider, however, the long exact homology sequence for the
pair $\,\left(M\x M,\,\Delta\right)\,$, which includes the
segment
\begin{equation*}
\cdots\ra H_n(\Delta)\ra H_n(M\x M)\ra
H_n(M\x M,\Delta)\ra\cdots\ .
\end{equation*}
By Lefschetz duality (cf. for example, \cite[Theorem 72.3]{mu}),
the last group above is isomorphic to $\,H^n\left((M\x M)\setminus\Delta\right)\,$,
which vanishes, as we just saw. The homomorphism
$\,H_n(\Delta)\ra H_n(M\x M)\,$ must then be surjective. But
this is impossible: $\,\Delta\approx M\,$ so
$\,H_n(\Delta)=\Z_2\,$, while $\,H_n(M\x M)\,$ has rank at
least two. This settles the case where $\,\Q\,$ is smooth, and
we have proven our Theorem.
\end{proof}

\bigskip


\vfill
\begin{thebibliography}{-1}

\parskip 4pt

\bibitem[G]{g} M.~Ghomi, {\it Non-existence of skew loops on
ellipsoids} To appear in PAMS. Preprint available at
{\tt http://front.math.ucdavis.edu/math.CO/0311292}$\ $.
\smallskip

\bibitem[GS]{gs} M.~Ghomi \& B.~Solomon,
{\it Skew loops and quadric surfaces}, Comment. Math. Helv.,
77 (2002) 767-782.
\smallskip

\bibitem[H]{h} M.~Hirsch, {\bf Differential Topology},
Springer-Verlag Graduate Texts in Mathematics {\bf 33}, New
York 1976.
\smallskip

\bibitem[Mi]{mi} J.~Milnor, {\bf Morse Theory}, Princeton
University Press 1963.
\smallskip

\bibitem[Mu]{mu} J.~Munkres, {\bf Elements of Algebraic Topology},
Addison-Wesley 1984.

\bibitem[Se]{se} B.~Segre, {\it Sulle coppie di tangenti fra
ioro parallele relative ad una curve chuisa sghemba}, {\bf
Hommage au Prof. Lucien Godeaux}, 141-167, Libraire
Universitaire, Louvain, 1968.
\smallskip

\bibitem[T]{t} S.~Tabachnikov, {\it On skew loops, skew branes
and quadratic hypersurfaces}, Moscow Math. J., 3 (2003),
681-690.
\smallskip

\bibitem[W]{w} J.H.~White, {\it Global Properties of
Immersions into Euclidean Spheres}, Indiana Univ. Math. J.,
{\bf 20} (12) 1971, 1187--1194
\smallskip

\end{thebibliography}
\end{document}